\documentclass{amsproc}
\usepackage{aspmproc, 
}

\newtheorem{thm}{Theorem}[section]
\newtheorem{lem}[thm]{Lemma}

\newtheorem{rem}[thm]{Remark}

\newcommand{\step}[1]{S{\small TEP}~#1}

\title[Correlations asymptotics
at high temperatures]{Random path representation and sharp 
correlations asymptotics at high-temperatures}

\author[M. Campanino]{Massimo Campanino
$^{1}$\footnote{$^{1}$Investigation supported by University of Bologna.                               Funds for selected research topics and by the MIUR                              national project "Stochastic processes ..."                                      }}
\author[D. Ioffe]{Dmitry Ioffe
$^{2}$\footnote{$^{2}$Supported by the Fund for the Promotion of Research 
at the Technion}}
\author[Y. Velenik]{Yvan Velenik}

\address{
Dipartimento di Matematica\\
Universit\`a di Bologna\\
Piazza di Porta S. Donato 5\\
I-40126  Bologna\\
Italy}
\address{
Faculty of Industrial Engineering\\
Technion\\
Haifa 3200\\
Israel
}
\address{
UMR-CNRS 6632\\
Universit\'e de Provence\\
39 rue Joliot Curie\\
13453 Marseille\\
France}

\rcvdate{XXXXXXX}% to be supplied by the volume editor
\rvsdate{XXXXXXX}% to be supplied by the volume editor

\begin{document}

\begin{abstract}
We recently introduced a robust approach to the derivation of sharp
asymptotic formula for correlation functions of statistical mechanics models in
the high-temperature regime. We describe its application to the nonperturbative
proof of Ornstein-Zernike asymptotics of 2-point functions for self-avoiding
walks, Bernoulli percolation and ferromagnetic Ising models. We then extend the
proof, in the Ising case, to arbitrary odd-odd correlation functions. We discuss
the fluctuations of connection paths (invariance principle), and relate the
variance of the limiting process to the geometry of the equidecay profiles.
Finally, we explain the relation between these results from Statistical
Mechanics and their counterparts in Quantum Field Theory.
\end{abstract}

\maketitle

\section{Introduction}
In many situations, various quantities of interest can be represented in terms
of path-like structures. This is the case, e.g., of correlations in various
lattice systems, either in perturbative regimes (through a suitable expansion),
or non-perturbatively, as in the ferromagnetic Ising models at supercritical
temperatures. Many important questions about the fine asymptotics of these
quantities can be reformulated as local limit theorems for these (essentially)
one-dimensional objects. In~\cite{CaIoVe02}, building upon the earlier
works~\cite{Io98,CaIo02}, we proposed a robust non-perturbative approach to such
a problem. It has already been applied successfully in the case of self-avoiding
walks, Bernoulli percolation and Ising models. We briefly review the results
that have been thus obtained (see also~\cite{CaIoVe03} for a short description of the main ideas of the proof).

\bigskip
\noindent\textbf{Self-avoiding walks.} A self-avoiding path $\omega$ from
$0$ to $x\neq 0$ is a sequence of distinct sites $t_0=0, t_1, t_2, \ldots,
t_n=x$ in $\mathbb{Z}^d$, with $|t_i-t_{i-1}|=1$, $i=1,\ldots n$ (the
restriction to nearest-neighbor jumps can be replaced by arbitrary, possibly
weighted, jumps of finite range). Let $\beta<0$, we are interested in the
following quantity:
$$
G_\beta^{\mathrm{SAW}}(x) \stackrel{\triangle}{=} \sum_{\omega:0\to x} e^{\beta |\omega|}\,,
$$
where the sum runs over all self-avoiding paths from $0$ to $x$, and $|\omega|$
denotes the length of the path. $G^{\mathrm{SAW}}_\beta(x)$ is finite for all
$\beta <\beta_{\mathrm{c}}^{\mathrm{SAW}}$, with
$\beta_\mathrm{c}^{\mathrm{SAW}}>-\infty$. Actually, $\sum_{x\in\mathbb{Z}^d}
G_\beta^{\mathrm{SAW}}(x)$ is finite if and only if
$\beta<\beta_{\mathrm{c}}^{\mathrm{SAW}}$.

\medskip
\noindent\textbf{Bernoulli bond percolation.} Let $\beta>0$. We consider a
family of i.i.d. $\{0,1\}$-valued random variables $n_e$, indexed by the
bonds $e$ between two nearest-neighbor sites of $\mathbb{Z}^d$ (again,
restriction to nearest-neighbor sites can be dropped);
$\mathrm{Prob}_\beta(n(e)=1) = 1-e^{-\beta}$. We say that $0$ is connected to
$x$ ($0 \leftrightarrow x$) in a realization $n$ of these random variables if
there is a self-avoiding path $\omega$ from $0$ to $x$ such that $n_e=1$ for all
increments $e$ along the path. We are interested in the following quantity: $$
G_\beta^{\textrm{perc}}(x) \stackrel{\triangle}{=} \mathrm{Prob}_\beta(0
\leftrightarrow x) \,.
$$
The high-temperature region $\beta<\beta_{\mathrm{c}}^{\textrm{perc}
}$ is
defined through 
\[
\beta_{\mathrm{c}}^{\textrm{perc}} \stackrel{\triangle}{=}
\sup\{\beta\,:\, \sum_{x\in\mathbb{Z}^d} G_\beta^{\textrm{perc} }(x) < \infty\}
> 0 .
\]
 It is a deep result of~\cite{AiBa87} that the  percolation transition 
is sharp, i.e.
\[
\beta_c\, =\, \inf\{ \beta \,:\, 
\mathrm{Prob}_\beta(0\leftrightarrow\infty)>0\}.
\]
\medskip
\noindent\textbf{Ising model.} Let $\beta>0$. We consider a
family of $\{-1,1\}$-valued random variables $\sigma_x$, indexed by the
sites $x\in\mathbb{Z}^d$. Let $\Lambda_L = \{-L,\ldots,L\}^d$. The probability
of a realization $\sigma$ of the random variables $(\sigma_x)_{x\in\Lambda_L}$,
with boundary condition $\overline{\sigma}\in\{-1,1\}^{\mathbb{Z}^d}$, is
given by
$$
\mu_{\beta,\overline{\sigma},L} (\sigma) \stackrel{\triangle}{=}
(Z_{\beta,\overline{\sigma},L})^{-1} \, \exp\bigl[ \beta \sum_{\substack{\{x,y\}
\subset \Lambda_L\\|x-y|=1}} \sigma_x\sigma_y + \beta
\sum_{\substack{x\in\Lambda_L,y\notin\Lambda_L\\|x-y|=1}}
\sigma_x\overline{\sigma}_y\bigr]\,.
$$
(As for the two previous models, the nearest-neighbor restriction can be
replaced by a -- possibly weighted -- finite-range assumption.) The set of
limiting measures, as $L\to\infty$ and for any boundary conditions, is a
simplex, whose extreme elements are the Gibbs states of the model. We define the
high-temperature region as $\beta<\beta_{\mathrm{c}}^{\textrm{Ising}}$, where
$$
\beta_{\mathrm{c}}^{\textrm{Ising}} = \sup \{\beta \,:\, \text{There is a
unique Gibbs state at parameter $\beta$}\} > 0\,.
$$
We are interested in the following quantity:
$$
G_\beta^{\textrm{Ising}}(x) \stackrel{\triangle}{=}
\mathbb{E}_{\mu_\beta}[\sigma_0 \sigma_x] \,,
$$
where the expectation is computed with respect to any translation invariant
Gibbs state $\mu_\beta$ (it is independent of which one is chosen). It is a deep
result of~\cite{AiBaFe87} that the high-temperature region can also be
characterized as the set of all $\beta$ such that
$$
\sum_{x\in\mathbb{Z}^d} G_\beta^{\textrm{Ising}}(x) < \infty\,.
$$

\bigskip
We now discuss simultaneously these three models; to that end, we simply forget
the model-specific superscripts, and simply write $\beta_{\mathrm{c}}$ or $G_\beta$.
It can be shown that for all three models, for all $\beta<\beta_{\mathrm{c}}$, the
function $G_\beta(x)$ is actually exponentially decreasing in $|x|$, i.e. the
corresponding inverse correlation length $\xi_\beta:\mathbb{R}^d\to\mathbb{R}$
satisfy 
$$
\xi_\beta(x) \stackrel{\triangle}{=} \lim_{k\to\infty} -\frac1{k} \,
\log G_\beta(\lfloor kx\rfloor ) > 0\,,
$$
where $\lfloor x \rfloor$ is the componentwise integer part of $x$.
Obviously, $\xi_\beta$ is positive-homogeneous, and it is not difficult to prove
that it is convex; it is thus an equivalent norm on $\mathbb{R}^d$ (for
$\beta<\beta_{\mathrm{c}}$).

The main result of~\cite{Io98, CaIo02, CaIoVe02} is the derivation of the
following sharp asymptotics for $G_\beta(x)$, as $|x|\to\infty$, for these three
models, in the corresponding high-temperature regions.
\begin{thm}
\label{thm_OZ}
Consider one of the models above, and let $\beta<\beta_{\mathrm{c}}$. 
Then, uniformly as $|x|\to\infty$,
$$
G_\beta(x) = \frac {\Psi_\beta(n_x)} {\sqrt{|x|^{d-1}}} \, e^{-\xi_\beta(n_x)\,
|x|}\, (1+o(1))\,,
$$
where $n_x = x/|x|$, and $\Psi_\beta$ is strictly positive and
analytic. Moreover, $\xi_\beta$ is also an analytic function.
\end{thm}
As a by-product of the proof of Theorem~\ref{thm_OZ}, we obtain the following
results on the shape of the equidecay profiles,
$$
\mathbf{U}_\beta \stackrel{\triangle}{=} \{ x\in\mathbb{R}^d \, :
\, \xi_\beta(x)
\leq 1 \}
$$
and their polar, the Wulff shapes
$$
\mathbf{K}_\beta \stackrel{\triangle}{=} \bigcap_{n\in\mathbb{S}^{d-1}} \{
t\in\mathbb{R}^d \,:\, (t,n)_d \leq \xi_\beta(n) \}\,.
$$
\begin{thm}
\label{thm_shapes}
Consider one of the models above, and let $\beta<\beta_{\mathrm{c}}$. 
Then $\mathbf{K}_\beta$ has a locally
analytic, strictly convex boundary. Moreover, the Gaussian curvature
$\kappa_\beta$ of $\mathbf{K}_\beta$ is uniformly positive,
\begin{equation}
\label{eq_PosStif}
\overline{\kappa}_\beta \stackrel{\triangle}{=}
\min_{t\in\partial\mathbf{K}_\beta} \kappa_\beta(t) > 0 \,.
\end{equation}
By duality, $\partial\mathbf{U}_\beta$ is also locally analytic and strictly
convex.
\end{thm}
\begin{rem}
In two dimensions $\mathbf{K}_\beta$ is reminiscent of the Wulff shape (and is
exactly the low-temperature Wulff shape in the cases of the nearest-neighbor
Ising and percolation models). Equation~\eqref{eq_PosStif} is then called the
positive stiffness condition; it is known to be equivalent to the following
sharp triangle inequality~\cite{Io94, PfVe99}: Uniformly in $u,v\in\mathbb{R}^2$
$$ \xi_\beta(u) + \xi_\beta(v) - \xi_\beta(u+v) \geq \overline{\kappa}_\beta
\left( |u| + |v| - |u+v| \right)\,.
$$
\end{rem}
\bigskip

Theorem~\ref{thm_OZ} can in fact easily be extended to arbitrary odd-odd correlation functions. We show this here in the most difficult case of ferromagnetic Ising models; namely, we establish exact asymptotic formula for correlation functions of the form
$
\mathbb{E}_{\mu_\beta} [ \sigma_A \sigma_{B+x} ]\,,
$
where $A,B$ are finite subsets of $\mathbb{Z}^d$ with $|A|$ and $|B|$ odd, and for any $C\subset \mathbb{Z}^d$, $\sigma_C \stackrel{\triangle}{=} \prod_{y\in C}\sigma_y$. Notice
that even-odd correlations are necessarily zero by symmetry. The case of even-even
correlations is substantially more delicate though (already for the much simpler SAW model), in particular in low dimensions; we hope to come back to this issue in the future.

\begin{thm}
\label{thm_oddodd}
Consider the Ising model. Let $\beta<\beta_{\mathrm{c}}^{\mathrm{Ising}}$, and
let $A$ and $B$ be finite odd subsets of ${\mathbb Z}^d$. 
 Then, uniformly in $|x|\to\infty$,
$$
\mathbb{E}_{\mu_\beta} [ \sigma_A \sigma_{B+x} ] =
\frac {\Psi_\beta^{A,B}(n_x)} {\sqrt{|x|^{d-1}}} \, e^{-\xi_\beta(n_x)\,
|x|}\, (1+o(1))\,,
$$
where $n_x = x/|x|$, and $\Psi_\beta^{A,B}$ is strictly positive and
analytic.
\end{thm}
We sketch the proof of this theorem in Section~\ref{sec_oddodd}.

\bigskip
The main feature shared by the three models discussed above is that the function $G_\beta(x)$ can each time be written in the form
\begin{equation}
\label{eq_1D}
G_\beta(x) = \sum_{\lambda:\, 0\to x} q_\beta (\lambda)\,,
\end{equation}
where the sum runs over admissible path-like objects (SAW paths, percolation clusters, random-lines, see Section~\ref{sec_RandLine}, respectively). The weights 
$q_\beta (\,\cdot\,)$ are supposed to be strictly positive and to possess a variation of the following four properties:
\begin{itemize}
\item \textbf{Strict exponential decay of the two-point function:} There 
exists $C_1<\infty$ such that, for
all $x\in\mathbb{Z}^d \setminus \{0\}$,
\begin{equation}
\label{eq:PUB}
g(x)\,  =\,  
\sum_{\lambda:\, 0\to x} q(\lambda) \leq C_1\, e^{-\xi(x)}\,,
\end{equation}
where $\xi(x) = -\lim_{k\to\infty}(k)^{-1}\, \log g(\lfloor kx\rfloor)$ is the inverse
correlation length.
\item\textbf{Finite energy condition:} For any pair of compatible 
paths $\lambda$ and $\eta$ define the conditional weight
$$
q(\lambda \,|\, \eta) = q(\lambda \amalg \eta ) / q(\eta)\,
$$
where $\lambda\amalg\eta$ denotes the concatenation of $\lambda$ and $\eta$. 
Then there exists a universal 
 finite constant $C_2<\infty$ such that the conditional
weights are controlled in terms of path sizes $|\lambda |$ as:
\begin{equation}
\label{eq:FE}
q(\lambda \,|\, \eta)\,  \geq \, e^{-C_2 |\lambda |}\, . 
\end{equation}
\item \textbf{BK-type splitting property:} There exists $C_3<\infty$, such 
that, for all
$x,y\in\mathbb{Z}^d \setminus \{0\}$ with $x\neq y$,
\begin{equation}
\label{eq:Spl}
\sum_{\lambda:\, 0\to x\to y} q(\lambda) \leq C_2\, \sum_{\lambda:\, 0\to x}
q(\lambda)\; \sum_{\lambda:\, x\to y} q(\lambda)\,.
\end{equation}
\item \textbf{Exponential mixing :}
There exists $C_4<\infty$ and $\theta\in (0,1)$ such that, 
for any four paths $\lambda$, $\eta$,
$\gamma_1$ and $\gamma_2$, with $\lambda \amalg \eta \amalg \gamma_1$ and
$\lambda \amalg \eta \amalg \gamma_2$ both admissible,
\begin{equation}
\label{eq:ExD}
\frac{q(\lambda \,|\, \eta\amalg\gamma_1)}{q(\lambda \,|\, \eta\amalg\gamma_2)}
\leq \exp\{C_4\, \sum_{\substack{x\in \lambda\\y\in\gamma_1\cup\gamma_2}}
\theta^{|x-y|} \}\,.
\end{equation}
\end{itemize}
Many other models enjoy a graphical representation of correlation functions of the
form~\eqref{eq_1D}. In perturbative regimes, cluster expansions
provide a generic  example. Non-perturbative examples include the random-cluster
representation for Potts (and other) models~\cite{FoKa72}, or
random walk representation of $N$-vector models~\cite{BrFrSp82}, etc...
However, it might not always be easy, or even possible, to establish
properties~\eqref{eq:PUB}, \eqref{eq:FE}, \eqref{eq:Spl} and~\eqref{eq:ExD} for the
corresponding weights, especially 
~\eqref{eq:Spl} which is probably the less robust one. It should 
however be possible to weaken the latter so that it 
only relies on some form of locally uniform mixing properties.

\subsection*{Road-map to the paper}
In Section~\ref{sec_fluct} we review and explain our probabilistic 
approach to the analysis of high temperature correlation functions. The 
point of departure is the random path representation formula  
\eqref{eq_1D}, and the whole theory is built upon a study of the 
local fluctuation structure of the corresponding connection paths. One of 
the consequences is the validity of the invariance principle under
 the diffusive scaling, which we formulate in Theorem~\ref{thm:invariance}
 below. For simplicity the discussion in Section~\ref{sec_fluct} is 
restricted to the case of SAW-s, and hence the underlying local limit
results are those about the sums of independent random variables. 
In the case of high temperature ferromagnetic Ising models the random 
line representation, which we shall briefly recall in 
Section~\ref{sec_RandLine}, gives rise to path weights $q_\beta$ which 
do not possess appropriate factorization properties. Nevertheless these 
weights satisfy conditions \eqref{eq:PUB}-\eqref{eq:ExD} and we conclude 
Section~\ref{sec_RandLine} with an explanation of how the problem 
of finding correlation asymptotics can be reformulated in terms
 of local limit properties of one dimensional systems generated by Ruelle 
operators for full shifts on countable alphabets. 
%In Section~\ref{sec_RandLine}, we 
%recall briefly the random-line representation
%of correlation functions in the ferromagnetic Ising models.
The proof of Theorem~\ref{thm_oddodd} is discussed in 
Section~\ref{sec_oddodd}. 
%In
%Section~\ref{sec_fluct}, we discuss the functional limit theorems which can be
%extracted from the Ornstein-Zernike asymptotics. 
Finally, in
Section~\ref{sec_QFT}, we explain the relation between the problems
discussed here, inspired by Statistical Physics, and their counterparts
originating from the corresponding lattice  Quantum Field Theories.

\section{Fluctuations of connection paths}
\label{sec_fluct}
In this section we describe local structure and large scale properties of 
connection paths conditioned to hit a distant point. In all three models above 
(SAW, percolation, Ising) the distribution of the connection  paths converges, after the 
appropriate rescaling, to the $(d-1)$-dimensional Brownian bridge, and, from the 
probabilistic point of view, these results belong to  the realm of classical Gaussian 
local limit analysis of one dimensional systems based on uniform analytic expansions
of finite volume log-moment generating functions.  An invariance principle for 
the sub-critical  Bernoulli bond percolation has been established in \cite{Ko}
 and for the phase separation line in the 2D nearest neighbour Ising model 
at any $\beta >\beta_c$ in \cite{GrIo}. In both cases the techniques and the ideas of \cite{CaIo02} and \cite{CaIoVe02} play the crucial role, and, in fact,
 the renormalization and the fluctuation analysis developed in the latter papers pertains
to a large class of models which admit a  random path type representation 
with path weights enjoying a suitable variation of \eqref{eq:PUB}\,-\eqref{eq:ExD}.
 In particular, it should lead to a closed form theory of low temperature phase 
boundaries in two dimensions \cite{IoVe}.
Note that different tools
have been early employed in \cite{Du, Hi}. 

For the sake of simplicity we shall sketch here the case of self-avoiding walks and
shall try to stipulate the impact of the geometry of ${\bf K}_\beta$ on the magnitude 
of paths fluctuations in the corresponding directions.

 Let $\hat{x}\in{\mathbb S}^{d-1}$ and the dual point $\hat{t}\in\partial {\bf K}_\beta ;\ 
(\hat{t}, \hat{x}) = \xi_\beta (\hat{x})$, 
 be fixed for the rest of the section. Consider the set ${\mathcal P}^n$ of all 
self-avoiding paths $\gamma :0\to \lfloor n\hat{x}\rfloor$, where 
 for   $y\in{\mathbb R}^d$  we define $\lfloor y \rfloor = 
(\lfloor y_1\rfloor ,\dots , \lfloor y_d\rfloor )\in {\mathbb Z}^d$. Finally, consider the 
following probability measure ${\mathbb P}^n_\beta$ on ${\mathcal P}^n$:
\begin{equation}
\label{eq:pm_saw}
{\mathbb P}^n_\beta  (\gamma )\, =\, 
\frac1{{\mathbf Z}_\beta^n } {\rm e}^{\beta  |\gamma |}  
1_{\{\gamma \in {\mathcal P}^n\}} .
\end{equation}
In order to explain and to formulate the invariance principle which holds under 
 ${\mathbb P}^n_\beta$ we need, first of all, to readjust the notion of 
irreducible splitting of paths $\gamma \in {\mathcal P}^n$;
\begin{equation}
\label{eq:gamma_split}
\gamma\, = \, \lambda_L\amalg \lambda_1\amalg\dots\amalg \lambda_M
\amalg\lambda_R .
\end{equation}
Fix $\delta \in (0,1)$ and a large enough renormalization scale $K$. Given 
a path $\lambda =(u_0, u_1, \dots ,u_m)$ let us say that a point 
$u_l ;\, 0<l< m$, is $\hat{x}$-correct break point of $\lambda$ if the following 
two conditions hold:
\medskip

{\bf A)} \quad $(u_j , \hat{x}) <(u_l , \hat{x}) <(u_i , \hat{x})$ for all $j< l< i$.

{\bf B)} \quad The remaining sub-path $(u_{l+1},\dots , u_m)$ lies inside the 
 set 
\[
2K {\mathbf U}_\beta (u_l ) + {\mathcal C}_\delta (\hat{t} ),
\]
where ${\mathbf U}_\beta (z) = z + {\mathbf U}_\beta$, and the forward cone ${\mathcal C}_\delta (\hat{t} )$ is defined as 
\begin{equation}
\label{eq:cone}
{\mathcal C}_\delta (\hat{t} )\, =\, \left\{ y\in {\mathbb R}^d\, :\, 
(y,\hat{t}) >(1-\delta )\xi_\beta (y)\right\}  .
\end{equation}
Note that this definition depends on the parameters $K$ and $\delta$; as they are usually kept constant, we only write them explicitly when needed. 
With $\hat{x}\in {\mathbb S}^{d-1}$, $\hat{t}\in\partial{\mathbf K}_\beta$,  
$K$ and $\delta$  fixed 
as above let us say that a path $\lambda$ is irreducible if it does not contain
 $\hat{x}$-correct break points.  We use ${\mathcal S}$ to denote the set 
of all irreducible paths (modulo ${\mathbb Z}^d$-shifts).  Define also the 
following three subsets of ${\mathcal S}$:
\begin{equation}
\label{eq:K_cone}
\begin{split}
 &{\mathcal S}_L = \{ \lambda = 
(u_0,\dots ,u_m)\in {\mathcal S}\, :\, \forall l>0\,
 (u_l ,\hat{x }) <(u_m ,\hat{x})\}\\[2ex]
&{\mathcal S}_R = \left\{
\begin{split}
 \lambda = (u_0,\dots ,u_m)&\in {\mathcal S}\, :\, \forall l>0\,
 (u_l ,\hat{x }) >(u_0 ,\hat{x})\ \  {\rm and}\\[2ex]
&\gamma\subset K{\mathbf U}_\beta (u_0 )
+{\mathcal C}_\delta (\hat{t} ) .
\end{split}
\right\}\\[2ex]
& {\mathcal S}_0 ={\mathcal S}_L\cap {\mathcal S}_R 
\end{split}
\end{equation}

For any $\gamma\in {\mathcal P}_n$ which has at least two $\hat{x}$-correct 
break points the 
 decomposition \eqref{eq:gamma_split} is unambiguously defined by
the following set of conditions:
\[
\lambda_L\in{\mathcal S}_L , \  \lambda_R\in{\mathcal S}_R\ {\rm and}\ 
\lambda_1,\dots ,\lambda_M\in {\mathcal S}_0 .
\]
The only difference between \eqref{eq:gamma_split}  and 
the irreducible decomposition employed in \cite{CaIoVe02} is that the 
break points here are defined with respect to the $\hat{x}$-orthogonal hyper-planes
 instead of $\hat{t}$-orthogonal hyper-planes. This is to ensure that the  displacements
along all the $\lambda$-paths which appear in \eqref{eq:gamma_split} have 
positive projection on the direction of $\hat{x}$. More precisely, given a SAW path
$\lambda = (u_0,\dots ,u_m)$ let us define the displacement along $\lambda$ 
as $V(\lambda ) = u_m -u_0$.  By the very definition of \eqref{eq:gamma_split}  all
\[
V_L\stackrel{\Delta}{=} V(\lambda_L ), \ V_1 
\stackrel{\Delta}{=} V(\lambda_1 ), \dots , 
\ V_M
\stackrel{\Delta}{=} V(\lambda_M ),\ 
V_R \stackrel{\Delta}{=} V(\lambda_R ) .
\]
belong to the (lattice) half-space $\{ y\in {\mathbb Z}^d\, :\, (y,\hat{x})>0\}$.
% As in the case of Lemma~\ref{lem:mass_gap_odd}  estimates of 
The renormalization calculus  developed in 
 \cite{CaIo02, CaIoVe02} implies:
\begin{lem}
\label{lem:mass_gap_path} For every $\beta <\beta_c$ and for any 
 $\delta >0$ there exists a finite scale $K_0 =
 K_0 (\delta , \beta )$ and a number $\nu =\nu (\delta , \beta )>0$, such that 
\begin{equation}
\label{eq:mass_gap_path}
\sum_{\lambda \in{\mathcal S}\, :\, V(\lambda) = y}\, {\rm e}^{\beta |\lambda |}
\ \leq \ {\rm exp}\left\{ -(\hat{t} ,y ) -\nu |y|\right\}\, ,
\end{equation}
uniformly in $y\in {\mathbb Z}^d$.
\end{lem}
Going back to the decomposition \eqref{eq:gamma_split} notice that 
\begin{equation}
\label{eq:sum_V}
V_L +V_1+\dots +V_M +V_R = \lfloor n\hat{x}\rfloor \, 
\end{equation}
for any 
$\gamma :0\to \lfloor n\hat{x}\rfloor$.  Therefore, Lemma~\ref{lem:mass_gap_path}
and the Ornstein-Zernike formula of Theorem~\ref{thm_OZ} yield:
\begin{equation}
\label{eq:V_control}
{\mathbb P}^n_\beta \left( \max\{|V_L| ,|V_1|,\dots, |V_M|,|V_R| \} >(\log n)^2\right)
\ =\ {\small o}\left( \frac1{n^\rho }\right) ,
\end{equation}
 for any $\rho >0$. In particular,  if for given $\gamma : 0\to 
\lfloor n\hat{x}\rfloor $
 one considers the  piece-wise constant  
trajectory $\widehat{\gamma}$ through the vertices  $0, V_L, V_L+V_1, \dots, 
\lfloor n\hat{x}\rfloor $, then the ${\mathbb R}^d$-Hausdorff distance between  
 $\gamma$ and $\widehat{\gamma}$ is bounded above as:
\begin{equation}
\label{eq:hausdorf}
{\mathbb P}^n_\beta \left( {\rm d}_{\rm H}
\left( \gamma ,\widehat{\gamma}\right)
> (\log n)^2\right)
\ =\ {\small o}\left( \frac1{n^\rho }\right) ,
\end{equation}
as well.  Indeed, one needs only to control the fluctuation of $\lambda_L$ in 
\eqref{eq:gamma_split}, the traversal deviations of paths in ${\mathcal S}_R$ are
automatically under control by the cone confinement property \eqref{eq:K_cone}.  
 
Estimate \eqref{eq:hausdorf} enables a formulation of the invariance principle 
for SAW $\gamma$ in terms of the effective path $\widehat{\gamma}$. In 
its turn the invariance principle for $\widehat{\gamma}$ is a version of the 
conditional invariance principle for paths of random walks in $(d-1)$-dimensions 
 with the direction of the target point $\hat{x}$ playing the role of time. It happens
 to be natural to choose the frame of the remaining $(d-1)$ spatial dimensions 
according to principal directions of curvature 
 ${\frak v}_1 , \dots ,{\frak v}_{d-1}$ of $\partial {\mathbf K}_\beta$ at $\hat{t}$.
In this way, in view of the positive $\hat{x}$-projection property of all the 
$\lambda$-path displacements  in \eqref{eq:gamma_split}, the effective path 
$\widehat{\gamma}\subset {\mathbb R}^d$ could be parametrized in 
the orthogonal frame $(\hat{x}, {\frak v}_1,\dots, {\frak v}_{d-1})$  as a function
$\widehat{X} :[0,n]\to {\mathbb R}^{d-1}$.  As usual define the diffusive 
scaling $\widehat{X}_n (\cdot )$ of $\widehat{X} (\cdot )$ as
\[
\widehat{X}_n (\tau )\, =\, \frac1{\sqrt{n}} \widehat{X} (\lfloor n\tau\rfloor ).
\]
\begin{thm}
\label{thm:invariance}
The distribution of $\widehat{X}_n (\cdot )$ under ${\mathbb P}^n_\beta$ weakly 
converges in $C_{0,0}[0,1]$ to the distribution of 
\begin{equation}
\label{eq:bridge}
\left( \sqrt{\kappa_1 } B_1 (\cdot ),\dots ,\sqrt{\kappa_{d-1} } B_{d-1} (\cdot )\right) ,
\end{equation}
where $B_1 (\cdot ),\dots , B_{d-1} (\cdot )$ are independent Brownian bridges 
on $[0,1]$ and $\kappa_1,\dots ,\kappa_{d-1}$ are the principal curvatures of 
$\partial {\mathbf K}_\beta$ at $\hat{t}$. 
\end{thm}
Let us dwell on the probabilistic picture behind Theorems ~
\ref{thm_OZ}, 
\ref{thm_shapes} and \ref{thm:invariance}:  
First of all, note that by Lemma~\ref{lem:mass_gap_path} 
\begin{equation}
\label{eq:Qnot}
{\mathbb Q}_0 ( y )\, =\, {\rm e}^{ (\hat{t} , y)} \sum_{\lambda \in 
{\mathcal S}_0\, :\, V(\lambda )=y} 
{\rm e}^{\beta |\lambda |}\ \stackrel{\Delta}{=} \ {\rm e}^{ (\hat{t} , y)}
{ W}_0 (y )
\end{equation}
is a (non-lattice) probability distribution on ${\mathbb Z}^d$ with exponentially 
decaying tails.  Indeed, an alternative
important way to think about  ${\mathbf K}_\beta$  is as of the closure of the 
domain of convergence of the series 
\begin{equation}
\label{eq:series}
t\in{\mathbb R}^d\, \mapsto \, \sum_{y\in {\mathbb Z}^d} {\rm e}^{(t,y )}G_\beta (y) .
\end{equation}
On the other hand, Lemma~\ref{lem:mass_gap_path} ensures that the series 
${\mathbb W}_0 (t) \stackrel{\Delta}{=} \sum {\rm e }^{(t,y)}W_0(y)$ converges in 
the $\nu$-neighbourhood $B_\nu  (\hat{t}) =\{ t\, :\, |t-\hat{t}|<\nu\}$ of $\hat{t}$.
In view of the decomposition \eqref{eq:gamma_split} and 
Lemma~\ref{lem:mass_gap_path}, 
\begin{equation}
\label{eq:renewal}
G_\beta ( n\hat{x}  )\, =\, O\left( {\rm e}^{-n\xi_\beta (\hat{x})
 -\nu n} \right)\ +\ \sum_{M=1}^{\infty} W_L\ast W_0^{\ast M}\ast  W_R 
(n\hat{x} ) ,
\end{equation}
where we have assumed for the convenience of notation that $n\hat{x}\in
{\mathbb Z}^d$, and 
$$
W_L (y) =\sum_{\lambda\in {\mathcal S}_L\, :\, V(\lambda)=y} 
{\rm e }^{\beta |\lambda |}\quad{\rm and}\quad
W_R (y) =\sum_{\lambda\in {\mathcal S}_R\, :\, V(\lambda)=y} 
{\rm e }^{\beta |\lambda |} .
$$
As a result,  the piece of the boundary  
$\partial {\mathbf K}_\beta$ inside $ B_\nu  (\hat{t})$ is implicitly given by 
\[
\partial {\mathbf K}_\beta\cap B_\nu  (\hat{t}) \, =\, 
\{ t\in B_\nu  (\hat{t}) \, :\, {\mathbb W}_0(t)=1\} .
\]
  In order to obtain the full claim of 
Theorem~\ref{thm_shapes} one needs only to check the non-degeneracy 
of ${\rm Hess}( {\mathbb W}_0)$ at $\hat{t}$, which, in the case of SAW-s, 
is a direct  consequence  of the finite energy condition 
\eqref{eq:FE}.    Note, by the way, that since 
$\hat{x}$ is the normal direction to $\partial {\mathbf K}_\beta$ at $\hat{t}$,
 there exists a number $\alpha \in (0,\infty )$, such that
\begin{equation}
\label{eq:xnormal}
\nabla {\mathbb W}_0 (\hat{t} )\, =\, \alpha \hat{x} .
\end{equation}
Multiplying both sides of \eqref{eq:renewal} by 
${\rm e}^{n\xi_\beta (\hat{x})} ={\rm e}^{(\hat{t} ,n\hat{x})}$ we arrive to the 
following key representation of the two point function $G_\beta$:
\begin{equation}
\label{eq:decomp}
\begin{split}
&{\rm e}^{n\xi_\beta (\hat{x})}G_\beta (n\hat{x}) \, 
=\, {\rm e}^{n(\hat{t}, \hat{x} )} 
 G_\beta ( n\hat{x})\ =\ {\small O}\left( {\rm e}^{-n\nu}\right)\ \\
%&+\, \sum_{M=1}^{\infty} \sum_{v_L ,v_R\in{\mathbb Z}^d} 
%W_L (y_L)W_R (y_R){\rm e}^{(\hat{t} , v_L+v_R )}
%{\mathbb Q}_0 (V_1+\dots +V_M =n\hat{x} - v_L -v_R )\\
&+\, \sum_{v_L ,v_R\in{\mathbb Z}^d} {\mathbb Q}_L (v_L){\mathbb Q}_R (v_R)
\sum_{M=1}^{\infty}  {\mathbb Q}_0 (V_1+\dots +V_M =n\hat{x} - v_L -v_R ) ,
\end{split}
\end{equation}
where, similar to \eqref{eq:Qnot}, we have defined ${\mathbb Q}_L (v) = 
{\rm e}^{(\hat{t} ,v)} W_L (v)$ and, accordingly, 
${\mathbb Q}_R (v) = 
{\rm e}^{(\hat{t} ,v)} W_R (v)$.

Unlike 
${\mathbb Q}_0$
the measures ${\mathbb Q}_L$ and ${\mathbb Q}_R$ are in general not
 probability but, by Lemma~\ref{lem:mass_gap_path}, 
 they are finite and have exponentially decaying tails:
\[
\sum_{|y| >n} \left({\mathbb Q}_L (y) + {\mathbb Q}_R (y)\right) \, 
\leq\, {\rm e}^{-\nu n/2} .
\]
Since by \eqref{eq:xnormal} the expectation of $V_l$ under ${\mathbb Q}_0$ 
equals to $\alpha \hat{x}$, the usual local limit CLT for ${\mathbb Z}^d$ 
random variables  and the Gaussian summation formula imply that the right
 hand side in \eqref{eq:decomp} equals to $c_1/\sqrt{n^{d-1}}$. Actually, a 
slightly more careful analysis along these line leads to the full analytic 
form of the Ornstein-Zernike formula as claimed in Theorem~\ref{thm_OZ}. 

Let us explain now how the principal curvatures $\kappa_1 ,\dots ,\kappa_{d-1}$
of $\partial {\mathbf K}_\beta$ at $\hat{t}$ enter the picture: By the 
irreducible path representation and arguments completely similar to those
just reproduced above,  the total weight of all piece-wise constant paths 
$\widehat{\gamma}=\widehat{\gamma}( V_L,V_1,\dots ,V_M ,V_R)\, ;\, M=1,2,\dots,$
which pass through a point $v_n\in {\mathbb Z}^d$;
\[
v_n \, =\, \lambda n \hat{x} \, +\, \sqrt{n} \sum_{l=1}^{d-1} a_l{\frak v}_l\, 
\stackrel{\Delta}{=}\,  \lambda n \hat{x} \, +\, \sqrt{n} {\frak v},
\]
equals to
\[ 
\frac{{c}_2}{\sqrt{\left( \lambda (1-\lambda ) n^2\right)^{(d-1)}}} 
{\rm e}^{-\xi_\beta (v_n) -\xi_\beta (n\hat{x} - v_n )} 
\left( 1\, +\, {\small o}(1)\right),
\]
where ${c}_2 >0 $ does not depend on $\lambda\in (0,1 )$ and 
 the coefficients $a_1,\dots ,a_{d-1}$.  Comparing with the OZ formula for the 
full partition function $G_\beta$, we infer that
\[
{\mathbb P}_\beta^n\left( v_n \in \widehat{\gamma}\right)\, =\, 
\frac{c_3\,  {\rm exp}\left\{
 -\left(\xi_\beta (v_n) +\xi_\beta (n\hat{x} - v_n )-\xi_\beta (n\hat{x})
\right)\right\}}{\sqrt{(\lambda  (1-\lambda ) n)^{d-1}}}
\left( 1\, +\, {\small o}(1)\right) .
\]
From now on we refer to Chapter~2.5 in \cite{Sch} for the missing details 
in the arguments below. $\xi_\beta$ is the support function of ${\mathbf K}_\beta$ 
 and by Theorem~\ref{thm_shapes} it is a smooth function. Thus, for every 
$v\in {\mathbb R}^d$ the gradient $\nabla \xi_\beta (v)\in \partial {\mathbf K}_\beta$
 and $\xi_\beta (v) = \left(  \nabla \xi_\beta (v) ,v\right)$ (in particular 
$\hat{t} =\nabla \xi_\beta (\hat{x}) =\nabla \xi_\beta (n\hat{x} )$).  Principal radii 
of the curvature $1/\kappa_1 ,\dots ,1/\kappa_{d-1}$ of $\partial {\mathbf K}_\beta$ at $\hat{t}$ are the eigenvalues of the linear map 
\[
{\rm d}^2\xi_\beta|_{\hat{x}} \,:\,T_{\hat{x}} {\mathbb S}^{d-1}\mapsto 
T_{\hat{x}} {\mathbb S}^{d-1} ,
\]
and ${\frak v}_1,\dots ,{\frak v}_{d-1}\in T_{\hat{x}} {\mathbb S}^{d-1}$ are the 
corresponding eigenvectors.  Therefore, 
\begin{equation*}
\begin{split}
\xi_\beta (v_n) &+\xi_\beta (n\hat{x} - v_n )-\xi_\beta (n\hat{x}) 
\,=\, n\lambda \left( \xi_\beta \left( \hat{x} +\frac1{\lambda\sqrt{n}} {\frak v}\right) - 
\xi_\beta (\hat{x} )\right)\,\\
& +\, 
n(1-\lambda )\left( \xi_\beta \left( \hat{x} -\frac1{(1-\lambda)\sqrt{n}} {\frak v}\right) - 
\xi_\beta (\hat{x} )\right)
\\
&=\, \frac1{2\lambda}\left( {\rm d}^2\xi_\beta|_{\hat{x}} {\frak v},{\frak v}\right)\, +
\, 
\frac1{2(1-\lambda)}\left( {\rm d}^2\xi_\beta|_{\hat{x}} {\frak v},{\frak v}\right)\, +
\, {\small O}\left( \frac1{\sqrt{n}}\right) \\
&=\, \frac1{2\lambda (1-\lambda )} \sum_{l=1}^{d-1} \frac{a_l^2 }{\kappa_l}
\, +
\, {\small O}\left( \frac1{\sqrt{n}}\right) .
\end{split}
\end{equation*}
Computations for higher order finite dimensional distributions follow a completely 
similar pattern.

\section{Random-line representation of Ising correlations}
\label{sec_RandLine}
Correlation functions of ferromagnetic Ising models admit a very useful
representation in terms of sums over weighted random paths, which is
especially  convenient for our purposes here. The two-point 
function formula \eqref{eq_1D} is a particular case.
 In this section, we recall 
how this representation is derived; we refer
to~\cite{PfVe99} for details and additional results. In the end of 
the section we shall briefly indicate how \eqref{eq:decomp} and, accordingly,
 the whole local limit analysis should be re-adjusted in order to incorporate  
the (dependent) case of Ising paths. 

\medskip
Although we use it for the infinite-volume Gibbs measure, it is convenient to
derive the random path 
 representation first for finite volumes, and then take the limit. As
there is a single Gibbs state for the values of $\beta$ we consider, it suffices
to consider free boundary conditions (i.e. no interactions between spins inside
the box and spins outside).

Given a set of edges $B$ of the lattice $\mathbb{Z}^d$, we define the associated
set of vertices as $V_B \stackrel{\triangle}{=} \{
x\in\mathbb{Z}^d \,:\, \exists\,e\in B\ \text{ with } x\in e\}$  ($x\in e$
means that $x$ is an endpoint of $e$). For any vertex $x\in V_B$, we define the
{\em index} of $x$ in $B$ by $\mathrm{ind}(x,B) \stackrel{\triangle}{=}
\sum_{e\in B} \mathbf{1}_{\{e \ni x\}}$. The {\em boundary} of $B$ is defined by
$\partial B \stackrel{\triangle}{=} \{x\in V_B \,:\, \mathrm{ind}(x,B) \text{ is
odd}\}$.

In this context, the finite volume Gibbs measure is defined by
$$
\mu_{B,\beta} (\sigma) \stackrel{\triangle}{=} Z_\beta(B)^{-1} \exp [ - \beta \sum_{e=(x,y)\in B} \sigma_x \sigma_y]\,,
$$
and we use the standard notation $\langle\,\cdot\,\rangle_{B,\beta}$ to denote expectation w.r.t. this probability measure.

We fix an arbitrary total ordering of $\mathbb{Z}^d$.
At each $x\in \mathbb{Z}^d$, we fix (in an arbitrary way) an
ordering  of the $x$-incident edges of the graph:
$$
B(x) \stackrel{\triangle}{=} \{e\in B \,:\, \mathrm{ind}(x,{\{e\}}) >0\} =
\{e^x_1\dots,e^x_{\mathrm{ind}(x,B)}\},
$$
and for two incident edges $e=e_i\in B(x)$, $e^\prime = e_j\in B(x)$ we say that
$e\leq e^\prime$ if the corresponding inequality holds for their sub-indices;
$i\leq j$.

\medskip
Let $A\subset V_B$ be such that $|A|$ is even; we write $\sigma_A
\stackrel{\triangle}{=} \prod_{i\in A} \sigma_i$. Using the identity $e^{\beta
\sigma_x \sigma_y} = \cosh(\beta) \bigl(1+\sigma_x\sigma_y \tanh(\beta)\bigr)$,
we obtain the following expression for the correlation function $\langle
\sigma_A \rangle_{B,\beta}$,
$$
\langle \sigma_A \rangle_{B,\beta} = Z_\beta(B)^{-1} \,
\sum_{\substack{D\subset B\\ \partial D = A}} \prod_{e\in D} \tanh
\beta\,,
$$
where
$$
Z_\beta(B) \stackrel{\triangle}{=} \sum_{\substack{D\subset B\\ \partial D = \emptyset}} \prod_{e\in D} \tanh\beta\,.
$$
From $D\subset B$ with $\partial D = A$, we would like to extract a family of
$|A|/2$ ``self-avoiding paths'' connecting pairs of sites of $A$. We apply the following algorithm:
\vskip 0.1cm
\noindent
\step{0} Set $k=1$ and $\Delta=\emptyset$.
\vskip 0.1cm
\noindent
\step{1}
Set $z^{(k)}_0$ to be the first site of $A$ in the ordering of ${\mathbb Z}^d$ fixed
 above,  $j=0$,
and update $A \stackrel{\triangle}{=} A\setminus \{z^{(k)}_0\}$.
\vskip 0.1cm
\noindent
\step{2}
 Let $e^{(k)}_j=(z^{(k)}_j,z^{(k)}_{j+1})$ be the first edge in
$B({z^{(k)}_j})\setminus\Delta $ (in the ordering of   $B({z^{(k)}_j})$
fixed above) such that $e^{(k)}_j\in D$. This defines $z^{(k)}_{j+1}$.
\vskip 0.1cm
\noindent
\step{3} Update $\Delta \stackrel{\triangle}{=} \Delta \cup \{e\in
B({z^{(k)}_j}) \,:\, e\leq e^{(k)}_j\}$.  If $z^{(k)}_{j+1}\in A$, then go to \step{4}.
Otherwise update $j \stackrel{\triangle}{=} j+1$ and return to \step{2}.
\vskip 0.1cm
\noindent
\step{4} Set $n^{(k)}=j+1$ and stop the construction of this path. Update
$A \stackrel{\triangle}{=} A\setminus \{z^{(k)}_{j+1}\}$, $k
\stackrel{\triangle}{=} k+1$ and go to \step{1}.

\vskip 0.1cm
This procedure produces a sequence
$(z^{(1)}_0,\ldots,z^{(1)}_{n^{(1)}},z^{(2)}_0,\ldots,z^{(|A|/2)}_{n^{(|A|/2)}})$.
Let $\overline i \stackrel{\triangle}{=} |A|/2+1-i$, and set $w^{(i)}_k \stackrel{\triangle}{=} z^{(\overline i)}_{n-k}$.

We, thus, constructed $|A|/2$ paths, $\gamma_i$, $i=1,\ldots,|A|/2$, given
by\footnote{This backward construction of the lines turns out to be convenient
for the reformulation in terms of Ruelle's formalism, see~\cite{CaIoVe02}.}
$$
\gamma_i \stackrel{\triangle}{=}
\gamma_i(D) \stackrel{\triangle}{=}
(w^{(i)}_0,\ldots,w^{(i)}_{n^{(i)}})
$$
connecting distinct pairs of points of $A$, and such that
\begin{itemize}
\item $(w^{(i)}_k,w^{(i)}_{k+1})\in B$, $k=0,\dots,n^{(i)}-1$,
$i=1,\ldots,|A|/2$
\item $(z^{(i)}_k,z^{(i)}_{k+1}) \neq (z^{(j)}_l,z^{(j)}_{l+1})$ if $i\neq j$,
or if $i=j$ but $k\neq l$.
 \end{itemize}
(but $z^{(i)}_k=z^{(j)}_l$ is allowed). 
A family of contours $\underline\gamma=(\gamma_1,\ldots,\gamma_{|A|/2})$ is $(A,B)$-admissible if it can be obtained from a set $D\subset B$ with $\partial D = A$, using this algorithm; in that case we write $\underline\gamma\sim (A,B)$. Notice that here the order of the paths is important: if $\gamma_k$ is a path from $x_k$ to $y_k$ then we must have $y_1 > y_2 > \ldots > y_{|A|/2}$. This is to ensure that we do not count twice the same configuration of paths.

The construction also yields a set of edges
$\Delta(\underline\gamma) \stackrel{\triangle}{=} \Delta$.
Observe that $\Delta(\underline\gamma)$ is entirely determined by $\underline\gamma$ (and the order chosen for the sites and edges). In particular the sets $D\subset B$ giving rise to an $(A,B)$-admissible family $\underline\gamma$ are characterized by $\partial D= A$ and
$$
D \cap \Delta(\underline\gamma) = \bigcup_{i=1}^{|A|/2} \gamma_i\,.
$$
Therefore, for such sets, $\partial (D\setminus \Delta(\underline\gamma)) = \emptyset$, and we can write
$$
\langle \sigma_A \rangle_{\beta,B} = \sum_{\underline\gamma \sim (A,B)} q_{\beta,B}(\underline\gamma)\,,
$$
where
$$
q_{\beta,B}(\underline\gamma) = w(\underline\gamma)\; \frac{Z_\beta(B\setminus \Delta(\underline\gamma))}{Z_\beta(B)}\,,
$$
with
$$
w(\underline\gamma) = \prod_{i=1}^{|A|/2} \prod_{k=1}^{n^{(i)}} \tanh\beta\,.
$$
This is an instance of the {\em random-line
representation} for correlation functions of the Ising model in $B$.
It has been studied in detail in \cite{PfVe97,PfVe99} and is essentially equivalent
(though the derivations are quite different) to the
random-walk representation of~\cite{Ai82}. We'll need a version of this
representation when $B$ is replaced by the set $\mathcal{E}(\mathbb{Z}^d)$ of all edges of $\mathbb{Z}^d$. To this end, we use the
following result (\cite{PfVe99}, Lemmas~6.3 and 6.9): For all $\beta<\beta_{\rm
c}$,
\begin{equation}
\langle \sigma_A \rangle_{\beta} = \sum_{\underline\gamma \sim A} q_\beta(\underline\gamma)\,,
\label{eq_RLTL}
\end{equation}
where $q_\beta(\underline\gamma) \stackrel{\triangle}{=} \lim_{B_n \nearrow \mathcal{E}(\mathbb{Z}^d)} q_{\beta,B_n}(\underline\gamma)$ is well
defined.

\medskip
It will also be useful to work with a more relaxed definition of admissibility, since we want to cut our paths into pieces, and the order of the resulting pieces might not correspond with the order of their endpoints. In general, given a path $\gamma=(x_1,x_2,\ldots x_n)$, we define $\Delta(\gamma) = \bigcup_{k=1}^n \{e\in B(x_k) \,:\, e\leq (x_{k-1},x_k) \}$. We say that a path $\gamma=(x_1,\ldots,x_n)$ is admissible if 
$\{ (x_1,x_2), \ldots, (x_{k-1}x_k)\}
 \cap \Delta((x_k,\ldots,x_n)) = \emptyset$ for all $2\leq k\leq n-1$.
Given a family of paths $\underline\gamma=(\gamma_1,\ldots,\gamma_n)$, we define $\Delta(\underline\gamma) = \bigcup_{k=1}^n \Delta(\gamma_k)$. A family of admissible paths $\underline\gamma$ is then admissible if 
$(\gamma_1,\ldots, \gamma_k) \cap \Delta((\gamma_{k+1},\ldots,\gamma_n)) = \emptyset$ for all $1\leq k\leq n-1$. Notice that the order of the paths is still important ($(\gamma_1,\gamma_2)$ can be admissible while $(\gamma_2,\gamma_1)$ is not), but there are no constraint on the order of their endpoints. Indeed, they can even share endpoints. Observe that these definitions are identical to those above when restricted to the same setting.

We then have the following crucial inequality: Let $\underline\gamma$ be an admissible family of paths. Then
\begin{equation}
\sum_{\substack{\gamma_0:x\to y\\ \gamma_0 \cap \Delta(\underline\gamma)=\emptyset}} q_\beta(\gamma_0,\underline\gamma) \leq q_\beta(\underline\gamma) \sum_{\gamma_0:x\to y} q_\beta(\gamma_0)\,.
\label{eq_BK}
\end{equation}
We give a brief proof. It is enough to consider the analogous statement in finite volumes $B$. Since $\Delta(\gamma_0,\underline\gamma) = \Delta(\gamma_0) \cup \Delta(\underline\gamma)$ and $\gamma_0 \cap \Delta(\underline\gamma) = \emptyset$, we have
$$
q_{\beta,B} (\gamma_0,\underline\gamma) = q_{\beta, B\setminus\Delta(\underline\gamma)} (\gamma_0) \, q_{\beta,B} (\underline\gamma)\,.
$$
Hence~\eqref{eq_BK} follows simply from Griffiths' second inequality since
$$
\sum_{\gamma_0 \subset B\setminus\Delta(\underline\gamma)} q_{\beta, B\setminus\Delta(\underline\gamma)} (\gamma_0) 
= \langle \sigma_x \sigma_y \rangle_{\beta, B\setminus \Delta(\underline\gamma)} 
\leq \langle \sigma_x \sigma_y \rangle_\beta\,.
$$
\medskip

If the set $A$ in \eqref{eq_RLTL} contains only two points, $A=\{ x,y\}$,
 then we recover \eqref{eq_1D}. 
The main difference between 
 the SAW case  considered in Section~\ref{sec_fluct}
 and the case of sub-critical 
ferromagnetic Ising models is that the path weights $q_\beta$ in  
\eqref{eq_1D} do not factorize: In general,
\[
 q_\beta (\gamma \amalg\lambda)\, \neq \, q_\beta (\gamma )
q_\beta (\lambda) .
\]
Consequently, 
  the 
displacement variables $V_1, V_2,\dots$ fail to be independent 
and   the underlying 
local limit analysis should be generalized.  
The appropriate framework is that of 
the statistical mechanics of one dimensional systems generated by 
Ruelle operators for full shifts on countable alphabets. We refer to 
\cite{CaIoVe02} for all the background material
 and here only sketch how the construction leads 
to the claims of Theorems~\ref{thm_OZ}, \ref{thm_shapes} and, after an
 appropriate re-definition of the measures ${\mathbb P}_\beta^n$, to the 
invariance principle stated in Theorem~\ref{thm:invariance}: As in the case 
of SAW-s fix a direction $\hat{x}\in {\mathbb S}^{d-1}$. The key 
renormalization result (Theorem~2.3 in \cite{CaIoVe02}) which implies that
 the rate of decay of the irreducible connections is strictly larger than the rate
of decay of the two point function $G_\beta$. In view of \eqref{eq:gamma_split}
 this validates a representation of $G_\beta$ as a sum of dependent 
random variables $V_L +V_1 +\cdots + V_M + V_R$ with exponentially 
decaying tails.  Namely, as in Section~\ref{sec_fluct} 
let ${\mathcal S}={\mathcal S}(K)$ be the set of all 
$\hat{x}$-irreducible paths. Then the following Ising analog of Lemma~\ref{lem:mass_gap_path} holds: 
\begin{lem}
\label{lem:mass_gap} For every $\beta <\beta_c$ and for any 
 $\delta >0$ there exists a finite  scale $K_0 =
 K_0 (\delta , \beta )$ and a number $\nu =\nu (\delta , \beta )>0$, such that 
\begin{equation}
\label{eq:mass_gap}
\sum_{\lambda \in{\mathcal S}\, :\, V(\lambda) = y}\,  q_\beta (\lambda )
\ \leq \ {\rm exp}\left\{ -(\hat{t} ,y ) -\nu |y|\right\}\, ,
\end{equation}
uniformly in $y\in {\mathbb Z}^d$.
\end{lem}
%Let, as before $\hat{x}=\vec{e}_1\in {\mathbb Z}^{d+1}$ be
% the direction of the imaginary time axis. 
The above Lemma suggests that the main contribution to the sharp 
asymptotics of $G_\beta$ comes from the weights of the paths 
$\lambda_1 ,\dots ,\lambda_M$ in the decomposition \eqref{eq:gamma_split}.
 Accordingly, consider now the set ${\mathcal S}_0$ of cylindrical  
$\hat{x}$-irreducible paths which was introduced in Section~\ref{sec_fluct}.
 Given a finite collection 
$\lambda ,\lambda_1 ,\dots ,\lambda_M \in {\mathcal S}_0$
 define the conditional weight
\[
q_\beta (\lambda\,  \big|\, \underline{\lambda})\, =\,
q_\beta (\lambda\,  \big|\, \lambda_1\amalg\dots \amalg \lambda_M )\, 
=\, \frac{
q_\beta (\lambda\amalg\lambda_1\amalg\dots \amalg \lambda_M )}
{q_\beta (\lambda_1\amalg\dots \amalg \lambda_M )} .
\]
By the crucial exponential mixing property \eqref{eq:ExD}
  one is able to control the dependence of the conditional weights 
$q_\beta (\lambda\,  \big|\, \underline{\lambda})$ on $\lambda_M$ as follows:
\begin{equation}
\label{eq:Lipshitz}
\sup_{\lambda ,\lambda_1,\dots ,\lambda_{M-1}\in {\mathcal S}_0}\,
\sup_{\lambda_M ,\tilde{\lambda}_M \in {\mathcal S}_0} 
\frac{q_\beta (\lambda\,  \big|\, \lambda_1\amalg\dots \amalg \lambda_M )}
{q_\beta (\lambda\,  \big|\, \lambda_1\amalg\dots \amalg \tilde{\lambda}_M )}\, 
\leq\, {\rm e}^{c_1 \theta^M} .
\end{equation}
In our formalism the 
set ${\mathcal S}_0$ plays the role of a countable alphabet. 
The estimate 
\eqref{eq:Lipshitz} enables the extension of the conditional weights
$q_\beta (\lambda\,  \big|\, \underline{\lambda})$  to the  case of infinite  
strings $\underline{\lambda} =(\lambda_1 ,\lambda_2 ,\dots )$. Let 
${\frak S}_{0, \theta}$ be the set of all such strings 
endowed with the metrics 
\[
{\rm d}_\theta (\underline{\lambda} ,\underline {\tilde{\lambda }})\, =\,
\theta^{\inf \{k :\lambda_k\neq \tilde{\lambda}_k\}} .
\]
and let ${\frak F}_{0, \theta}$ be the set 
of all bounded Lipschitz continuous functions 
on ${\frak S}_{0, \theta}$. 

 As before we choose  
 $\hat{t}
% = (\hat{p_1},0)
\in \partial {\mathbf K}_\beta$
  to be the dual direction to
% the imaginary time axis 
$\hat{x}$.  
  Given a path $\lambda\in {\mathcal S}_0$ and a string 
$\underline{\lambda}\in {\frak S}_{0, \theta}$ define the potential
\[
\psi_\beta (\lambda\, \big|\, \underline{\lambda})\, =\, 
\log q_\beta (\lambda\, \big|\, \underline{\lambda}) + 
\left(\hat{t} ,V(\lambda )\right) .
\]
By \eqref{eq:ExD} and  Lemma~\ref{lem:mass_gap} the operator 
\begin{equation}
\label{eq:Ruelle}
{\mathcal L}_{z} f(\underline{\lambda })\, =\, 
\sum_{\lambda\in {\mathcal S}_0} 
{\rm e}^{ \psi_\beta (\lambda\, |\, \underline{\lambda}) + (z, V(\lambda ))}
 f(\lambda\amalg \underline{\lambda }) , 
\end{equation}
is well defined and bounded on ${\frak F}_{0,\theta} $ for  
every 
 $z\in {\mathbb C}^d $ with $  |z|<\nu$. 
\smallskip

The dependent Ising analog of \eqref{eq:decomp} 
is then given (see Section~3 of \cite{CaIoVe02}) by
\begin{equation}
\label{eq:Gdecomp}
\begin{split}
&{\rm e}^{n\xi\beta (\hat{x})}G_\beta (n\hat{x}) \
 =\ {\small O}\left( {\rm e}^{-n\nu}\right)\ \\
%&+\, \sum_{M=1}^{\infty} \sum_{v_L ,v_R\in{\mathbb Z}^d} 
%W_L (y_L)W_R (y_R){\rm e}^{(\hat{t} , v_L+v_R )}
%{\mathbb Q}_0 (V_1+\dots +V_M =n\hat{x} - v_L -v_R )\\
&+\, 
\sum_{\mu\in {\mathcal S}_L}\sum_{\eta\in {\mathcal S}_R} q_\beta (\mu )
q_\beta (\eta ) 
%\sum_{v_L ,v_R\in{\mathbb Z}^d} {\mathbb Q}_L (v_L){\mathbb Q}_R (v_R)
\sum_{M=1}^{\infty} 
 {\mathbb Q}_{0,M}^{\mu ,\eta} (n\hat{x} - v_L -v_R ) .
\end{split}
\end{equation}
For each $M=1,2,\dots$ the family of weights 
$\{{\mathbb Q}_{0,M}^{\mu ,\eta}\}$ is related to 
the family of operators $\{ {\mathcal L}_z\}$ via the Fourier transform:
\begin{equation}
\label{eq:Rdecomp}
\sum_{y\in {\mathbb Z}^d} {\rm e }^{(z ,y)} 
{\mathbb Q}_{0,M}^{\mu ,\eta} (y) \, =\, 
{\mathcal L}_{z}^M w_{\mu , \eta} ,
\end{equation}
where the family $\{ w_{\mu ,\eta}\}$ is uniformly positive and uniformly 
bounded in ${\frak F}_{0,\theta}$. In this way 
 the analytic perturbation theory 
of the leading (that is lying on the spectral circle) eigenvalue of 
${\mathcal L}_z$ enables the expansion of the logarithm of the right hand
side in \eqref{eq:Gdecomp} which, in its turn, leads to classical 
Gaussian local
limit results for the dependent sums $V_1 +\cdots + V_M$.

\section{Asymptotics of odd-odd correlations}
\label{sec_oddodd}
In this section, we sketch the proof of Theorem~\ref{thm_oddodd}. We do not give a complete,
self-contained argument, since this would be too long, and would involve many repetitions from~\cite{CaIoVe02}. Instead, we
provide the only required update as compared to the proof for $2$-point
functions given in the latter work. As such, this section should be considered
as a complement,
%. We use the same notations, and do not
%redefine here all the needed quantities, see~
 and we shall give exact references to the formulas in \cite{CaIoVe02} whenever 
required.

As explained in Section~\ref{sec_RandLine}, the correlation function $\langle
\sigma_A \sigma_{B+x} \rangle_\beta$ admits a random-line representation of the
form
$$
\langle \sigma_A \sigma_{B+x}
\rangle_\beta = \sum_{\underline{\gamma} \sim A\cup (B+x)}
q_\beta(\underline{\gamma})\,,
$$
where $\underline\gamma$ runs over families of
compatible open contours connecting all the sites of $A\cup (B+x)$.
Among the
$\tfrac12(|A|+|B|)$ paths of $\underline{\gamma}$, at least one must connect a
site of $A$ to a site of $B+x$. We first show that one can ignore the
contribution of $\underline{\gamma}$ with more than one such connection (i.e. at
least three of them). The first observation is that we have the following lower
bound on the correlation function: By the second Griffiths' inequality,
$$
\langle \sigma_A \sigma_{B+x} \rangle_\beta \geq \langle \sigma_{A\setminus
\{y\}}\rangle_\beta \, \langle\sigma_{B\setminus\{z\}} \rangle_\beta \, \langle
\sigma_y \sigma_{z+x} \rangle_\beta\,,
$$
where $y$ and $z$ are arbitrarily chosen sites of $A$ and $B$ respectively.
Another application of the second Griffiths' inequality implies that
$$
\langle \sigma_{A\setminus \{y\}}\rangle_\beta\, \langle
\sigma_{B\setminus\{z\}} \rangle_\beta > 0\,.
$$
Moreover, we already know that
$$
\langle \sigma_y \sigma_{z+x} \rangle_\beta = \Psi_\beta
(n_x)\, |x|^{-(d-1)/2}\, e^{-\xi_\beta(x)}\; (1+o(1))\,.
$$
But, applying~\eqref{eq_BK},
we obtain immediately that the contribution of families of paths
$\underline{\gamma}$ with three or more connections between $A$ and $B+x$ is
bounded above by $C(A,B)\, e^{-3\xi_\beta(x)}$ and is therefore negligible.

\smallskip
We can henceforth safely assume that there is a single connection between $A$
and $B+x$; we denote the corresponding path by $\gamma$, while the remaining
paths are denoted by $\underline\gamma_A$ and $\underline\gamma_B$. We want to
show that we can repeat the argument used for the two-point
function in~\cite{CaIoVe02} in this more general setting. This is indeed quite reasonable since the paths in $\underline\gamma_A$ and $\underline\gamma_B$ should remain localized, and therefore the picture is still that of a single very long path as for 2-point functions. The main point is thus to prove sufficiently strong localization properties for the paths $\underline\gamma_A$ and $\underline\gamma_B$, so as to ensure that an
 appropriate  version Lemma~\ref{lem:mass_gap}  (see also  Theorem~2.3
 of~\cite{CaIoVe02}) still holds. The import of the latter lemma was to assert 
 nice decay and decoupling properties of the integrated 
 weights of the  irreducible pieces in the decomposition
 of connection  paths \eqref{eq:gamma_split}.
%the construction of a decomposition of the random path into a collection of irreducible pieces with .
 Notice first that exactly the same decomposition can still be used here, provided we attach the paths in $\underline\gamma_A$ and $\underline\gamma_B$ to the corresponding leftmost and rightmost 
extremal pieces $\lambda_L$ and $\lambda_R$,  and keep the remaining 
intermediate cylindrical irreducible
pieces  unchanged. Apart from the compatibility requirements one then has to check that $\underline\gamma_B$ stays
inside the forward cone containing $\lambda_R$, so that the crucial estimate~(3.9) in~\cite{CaIoVe02} remains valid.

\bigskip
Let $y\in A$, $z\in B+x$, $\gamma:\, y\to z$, and let $\underline\gamma_A$, resp. $\underline\gamma_B$, denote the collections of 
remaining paths connecting pairs of sites in $A\setminus\{y\}$, respectively $x+B\setminus\{z\}$. For given collections 
$\underline\gamma_A$ and $\underline\gamma_B$ we 
 define the irreducible decomposition of $\gamma$ in precisely the same way as in~
\eqref{eq:gamma_split}, except for the extremal pieces 
$\lambda_L = (u^L_0 ,\dots ,u^L_m )$ and 
$\lambda_R = (u^R_0 ,\dots ,u^R_n )$, 
which have to satisfy the following modified set of conditions:
\begin{itemize}
\item $(u^L_k ,\hat{x}) < (u^L_m,\hat{x}) \ \forall\, k=0,\dots ,m-1$
\item $(u^R_k ,\hat{x}) > (u^R_0,\hat{x}) \ \forall\, k=1,\dots , n$
\item $\underline\gamma_A$ must belong to the
 same $\hat{x}$-halfspace as $\lambda_L$ and for any $\hat{x}$-break point $u^L_k$ of
$\lambda_L$ the $\hat{x}$-orthogonal hyperplane through $u^L_k$ intersects 
$\underline\gamma_A$.
\item $\underline\gamma_B$ must belong to the
 same $\hat{x}$-halfspace as $\lambda_R$ and for any $\hat{x}$-break point $u^R_k$ of
$\lambda_R$ the $\hat{x}$-orthogonal hyperplane through $u^R_k$ intersects 
$\underline\gamma_B$.
\item $\underline\gamma_B$ must belong to $2K\mathbf{U}_\beta( u^R_0) + \mathcal{C}_\delta({t})$ (see 
  \eqref{eq:cone}).
\end{itemize}
With a slight ambiguity  of notation let us call compatible pairs 
$(\underline\gamma_A ,\lambda_L )$ and  $(\underline\gamma_B ,\lambda_R )$
 $\hat{x}$-irreducible if they satisfy all the conditions above. 
We then only have to check that
\begin{gather*}
\sum_{\substack{(\underline\gamma_A,\lambda_L)\,\text{$\hat{x}$-irreducible}\\
\lambda_L :y\to u}} q_\beta(\underline\gamma_A,\lambda) \leq e^{-(\hat{t},u-y) -\nu |u-y|}\,, \\
\sum_{\substack{(\underline\gamma_B,\lambda_R)\,\text{$\hat{x}$-irreducible}\\
\lambda_R :u\to z}} q_\beta(\underline\gamma_B,\lambda) \leq e^{-(\hat{t},z-u) -
\nu |z-u|} 
\end{gather*}
for some $\nu>0$ and any $u\in\mathbb{Z}^d$. We only 
check the second statement since it is the most complicated one.
 Fix a large enough scale $K$.  A site $u$
of $\gamma$ is a $(\hat{x}, \gamma_B ,\delta)$-admissible break point if it
is a $\hat{x}$-break point of $\gamma$ and, in addition, 
$$
\underline\gamma_B 
\subset 2K\mathbf{U}_\beta(u) + \mathcal{C}_\delta(\hat{t})\,.
$$
\begin{lem}
\label{lem:mass_gap_odd}
Fix a forward cone parameter $\delta \in (0,1/4 )$ and a set
$B = \{y_1,z_1, \ldots, y_n, z_n\} ;\ 
B\subset\mathbb{Z}^d\setminus \{0\}$. There exist a renormalization scale
$K_0$ and positive numbers $\epsilon =\epsilon (\delta ,\beta )$, $\nu =\nu
(\delta ,\beta )$ and $N=N(\beta)<\infty$,
such that for all $K\geq K_0$, the upper bound
 \begin{equation*}
\sum_{\substack{\lambda :-x\to 0\\
\partial\underline\gamma_B}} q_\beta(\lambda  ,\underline\gamma_B )
\mathbf{1}_{ \{ \text{$\lambda $ has no $(\hat{x},\underline\gamma_B)$-admissible
break points} \} } \leq  N \mathrm{e}^{-(t,x)_d -\nu |x|} \,,
\end{equation*}
holds uniformly in the dual directions $t\in\partial\mathbf{K}_\beta$ and in the
starting points $x\in\mathbb{Z}^d$. In the first sum $\underline\gamma_B,\lambda$ runs over all admissible family of paths such that $\lambda:\, -x\to 0$, while $\underline\gamma_B = (\gamma_1,\ldots,\gamma_n)$ satisfies $\gamma_k:\, y_k\to z_k$.
\end{lem}
\begin{proof}
Applying~\eqref{eq_BK}, we can assume that
$x\in\mathcal{C}_\nu'(t)$, see the remark after Theorem~2.3 of~\cite{CaIoVe02}. Let $Q_C(x) =
\{v\in\mathbb{Z}^d \,:\, |v|\leq |x|/C\}$ where $C$ is some large enough constant.
To simplify notations, we suppose that $B=\{y,z\}$, i.e. that $\underline\gamma_B
\equiv \gamma:\, y\to z$. The general case is treated in the same way.

We first show that, typically, $\gamma \subset Q_C(x)$. Indeed, using
again~\eqref{eq_BK}, we have that
\begin{align*}
\sum_{\substack{(\lambda,\gamma)\\
\lambda:\, -x\to 0\\
\gamma:\, y\to z }} 
q_\beta(\lambda,\gamma)
\mathbf{1}_{ \{ \gamma \not\subset Q_C(x)  \} }
&\leq
\sum_{u\in\partial Q_C(x)}
\sum_{\substack{(\lambda,\gamma_1,\gamma_2)\\
\lambda:\, -x\to 0\\
\gamma_1:\, y\to u,\, \gamma_2:u\to z }} 
q_\beta(\lambda,\gamma_1,\gamma)\\
&\leq
\sum_{u\in\partial Q_C(x)}
\langle \sigma_{-x} \sigma_0 \rangle_\beta\;
\langle \sigma_y \sigma_u \rangle_\beta\;
\langle \sigma_u \sigma_v \rangle_\beta\\
&\leq
\frac{c_d|x|}{C}\, \mathrm{e}^{-c|x|}\, e^{-\xi_\beta(x)}\,.
\end{align*}
We can therefore suppose that $\gamma\subset Q_C(x)$.
Observe now that in the latter case
\begin{multline*}
\{\lambda \text{ has no } ({x}, \gamma, 2\delta) \text{-break
points} \}\\
\subset \{\lambda \text{ has no } {x}
\text{-break point $u$ with } ({t},u)  \leq -\tfrac12 ({t},x)
\}\\
\stackrel{\triangle}{=} \mathcal{A}(t,K,\delta,x)\,,
\end{multline*}
provided that $C$ is
taken large enough. Indeed, would such a ${x}$-break point
$u$ exist then the 
cone $u+\mathcal{C}_{2\delta}(t)$ must contain the box $Q_C(x)$,
hence also $\gamma$.

The probability of $\mathcal{A}(t,K,\delta,x)$ is estimated exactly as in the proof of Theorem~2.3 of~\cite{CaIoVe02}. Indeed, the presence of the path $\gamma$ only affects an arbitrarily small fraction of the slabs 
$\mathcal{S}_k(t)$ introduced in the latter proof, provided 
$C$ is taken large enough, so that the argument given there 
applies with no modifications.
\end{proof}

\section{Relation to Quantum Field Theories}
\label{sec_QFT}
 There is an abundant literature devoted to 
the relation between Ising and other ferromagnetic type models to the Euclidean lattice
quantum field theories, see e.g. \cite{Schor, P-L} or more recently \cite{BarF, AuB}; 
the latter article contains also an extensive bibliography on the subject. In this works the spins live on the integer lattice ${\mathbb Z}^{d+1}$  with one special direction, 
say $\vec{e}_1$, 
 being visualized as the imaginary time axis. Thus, for example, the analyticity properties of the mixed Fourier transform 
\begin{equation}
\label{eq:Fourier}
{\mathbb G}_\beta (p_1, i{\mathbf p})\, =\, \sum_{x_1\in {\mathbb Z}}
\sum_{{\mathbf x}\in{\mathbb Z}^d} {\rm e}^{p_1x_1 + i({\mathbf p},{\mathbf x})}
G_\beta (x_1 ,{\mathbf x}) ,
\end{equation}  
$(p_1 ,{\mathbf p})\in {\mathbb T}\times {\mathbb T}^d$, are related in this way to the 
question of existence of one particle states.

Below we shall briefly indicate how the the key probabilistic representation 
\eqref{eq:decomp} leads to the following conclusion (see e.g Proposition~4.2 in
 \cite{P-L}, Theorem~2.3 in \cite{Schor}): For every ${\mathbf p}\in 
{\mathbb T}^d$ define 
\begin{equation}
\label{eq:mass_shell}
\omega ({\mathbf p})\, =\, -\lim_{n\to\infty}\frac1n \log\sum_{{\mathbf x\in {\mathbb Z}^d}}
{\rm e}^{i ({\mathbf p},{\mathbf x})} G_\beta (n,{\mathbf x}) ,
\end{equation}
 $\omega ({\mathbf p})$ being interpreted as the energy of a particle with momentum 
${\mathbf p}$.
\begin{thm}
\label{thm:omega_p}
There exists a neighbourhood $B_\delta =\{ {\mathbf p}\, :\, |{\mathbf p}|<\delta\}$ 
of the origin in ${\mathbb R}^d$ such that the function ${\mathbf p}\mapsto 
\omega ({\mathbf p})$ is real analytic on $B_\delta$. 
${\rm Hess}\left(\omega\right) (0)$ is precisely the matrix of the 
 second fundamental form of 
$\partial {\mathbf K}_\beta$ at $\hat{t} = (\xi_\beta (\vec{e}_1), 0)$. 
 Furthermore, there exists $\epsilon >0$ 
% a neighbourhood 
%$B_{\epsilon }^{{\mathbb C}} (\hat{p}_1) =\{ P\, :\, |p-\hat{p}_1 |<\epsilon \}$ of  
% $\hat{p}_1 =  \xi_\beta (\vec{e}_1 )$ in ${\mathbb C}$ 
such that 
for every ${\mathbf p}\in B_\delta$ the function 
\[
p_1\, \mapsto \, {\mathbb G}_\beta (p_1, i{\mathbf p})
\]
has a meromorphic extension to the disc 
$\{ p_1\in {\mathbb C}\, :\, |p_1 -\hat{p}_1 |<\epsilon\} ;\,  \hat{p}_1 =  \xi_\beta (\vec{e}_1 )$, with the 
 only  simple pole at $p_1 =\omega ({\mathbf p})$. 
\end{thm}
In the sequel we use the notation introduced in Section~\ref{sec_fluct}. Because of the 
${\mathbb Z}^d$-lattice symmetries the dual point $\hat{t}\in \partial {\mathbf K}_\beta$
of $\vec{e}_1$ is given by $\hat{t}= (\hat{p}_1 ,0)$ with  
$\hat{p}_1 =\xi_\beta (\vec{e}_1 )$. 
 Given $y\in {\mathbb Z}^{d+1}$ define (see \eqref{eq:K_cone})
\[
W(y)\, =\, \sum_{\substack{\gamma :0\to y\\ \gamma\in {\mathcal S}}}
{\rm e}^{\beta |\gamma |} ,
\]
and let $W_L ,W_0$ and $W_R$ be defined as in Section~\ref{sec_fluct}.
Summing up all the weights of irreducible paths in \eqref{eq:gamma_split} 
we arrive to the following representation of $G_\beta $:
\begin{equation}
\label{eq:G_representation}
\begin{split}
G_\beta (y)\, &=\, W(y) + \sum_{y_L +y_R =y} W_L (y_L)W_R(y_R) \\
&+\, 
\sum_{y_L ,y_R}\sum_{M=1}^\infty W_L(y_L)W_R(y_R) W^{*M}_0 (y-y_L -y_R) .
\end{split}
\end{equation}
Consider the mixed Fourier transforms
\begin{equation*}
\begin{split}
&{\mathbb W}(p_1 ,{\mathbf p})\, =\, \sum_{x_1\in{\mathbb Z}}
\sum_{{\mathbf x}\in {\mathbb Z}^d}
 {\rm e}^{p_1 x_1 +({\mathbf p},{\mathbf x})} W(x_1, {\mathbf x})\\
&\qquad {\rm and}  \\
&{\mathbb W}_b(p_1 ,{\mathbf p})\, =\, \sum_{x_1\in{\mathbb Z}}\sum_{{\mathbf x}\in 
{\mathbb Z}^d} {\rm e}^{p_1 x_1 +({\mathbf p},{\mathbf x})} W_b(x_1, {\mathbf x})\, ;\ 
b=0,L,R .
\end{split}
\end{equation*}
By Lemma~\ref{lem:mass_gap_path} all four functions above are analytic in the 
complex neighbourhood $B_\nu^{{\mathbb C}} (\hat{t} )$ of $\hat{t}$; 
$B_\nu^{{\mathbb C}}  (\hat{t} ) =\{ (p_1 ,{\mathbf p}):\sqrt{|p_1 - \hat{p}_1|^2
 +|{\mathbf p}|^2} <\nu \}$. Thus,
the extension of ${\mathbb G}_\beta (p_1 ,{\mathbf p})$ 
to  $B_\nu^{{\mathbb C}} (\hat{t} )$ is 
given by:
\[
{\mathbb W}(p_1 ,{\mathbf p})\, +\, \frac{  {\mathbb W}_L (p_1 ,{\mathbf p})
{\mathbb W}_R(p_1 ,{\mathbf p})}
{1- {\mathbb W}_0(p_1 ,{\mathbf p})}.
\]
Consequently, the surface of poles of ${\mathbb G}_\beta (p_1 ,{\mathbf p})$  
inside $B_\nu^{{\mathbb C}} (\hat{t} )$ is given by the implicit equation 
\begin{equation}
\label{eq:poles}
{\mathbb W}_0(p_1 ,{\mathbf p})=1. 
\end{equation}
As we have already seen  in  Section~\ref{sec_fluct}, the restriction of \eqref{eq:poles}
 to $(p_1 ,{\mathbf p})\in {\mathbb R}\times {\mathbb R}^d$ defines the piece of 
the boundary $\partial {\mathbf K}_\beta$ inside $B_\nu (\hat{t} )$.  Since by 
\eqref{eq:xnormal} $\partial {\mathbb W}_0 /\partial p_1 (\hat{t})\neq 0$ and, in 
addition, ${\rm Hess}({\mathbb W}_0) (\hat{t})$ is non-degenerate, the analytic 
implicit function theorem implies that there exists $\delta > 0$, such that the 
equation \eqref{eq:poles} can be resolved for 
${\mathbf p}\in B_\delta^{{\mathbb C}}(\hat{t})\subset {\mathbb C}^d$ as 
\begin{equation}
\label{eq:poles_resolved}
p_1\, =\, \widetilde{\omega }({\mathbf p} ).
\end{equation}
In particular, for ${\mathbf p}\in {\mathbb R}^d$ the equation 
\eqref{eq:poles_resolved} gives a parameterization of $\partial {\mathbb K}_\beta$
 in the $\delta$-neighbourhood of $\hat{t}$, and ${\rm Hess}(\widetilde{\omega})(0)$
is, indeed, the matrix of the second fundamental form of $\partial  {\mathbf K}_\beta$
 at $\hat{t}$. Finally, the  1-particle mass shell $\omega$ in 
 \eqref {eq:mass_shell} is recovered as
 $\omega ({\mathbf p}) =\widetilde{\omega }(i {\mathbf p})$.

\medskip

In the case of ferromagnetic Ising models set $\vec{p} = (p_1, {\mathbf p})$ and 
readjust the definition \eqref{eq:Ruelle} of the Ruelle operator ${\mathcal L}_z$ as
$
\widetilde{{\mathcal L}}_{\vec{p}} = {\mathcal L}_{\vec{p}-\hat{t}}$. 
Then $\widetilde{{\mathcal L}}_{\vec{p}}$ 
is well defined and bounded on ${\frak F}_{0,\theta} $ for  every 
 $\vec{p} \in B_\nu^{{\mathbb C}} (\hat{t} )$. It could be then shown 
that the surface of poles of ${\mathbb G}_\beta$ inside 
$B_\nu^{{\mathbb C}} (\hat{t} )$ is implicitly given by
\[
\tilde{\rho}_{\beta }( p_1, {\mathbf p})\, =\, 1 , 
\]
 where $\tilde{\rho}_\beta (p_1 ,{\mathbf p} )$ is the leading  (lying on the 
spectral circle) eigenvalue of $\widetilde{{\mathcal L}}_{\vec{p}}$. Further analysis 
of the spectral properties of the family $\{\widetilde{ {\mathcal L}}_{\vec{p}}\}$ reveals 
\cite{CaIoVe02} that there exists $\epsilon >0$, such that $\rho_\beta (p_1 ,{\mathbf p} )$ is a simple pole of the 
corresponding resolvent  for every $\vec{p}\in B_\epsilon^{{\mathbb C}} (\hat{t} )$. 
 In this way, the conclusion of Theorem~\ref{thm:omega_p} follows 
from the analytic 
 perturbation theory of discrete spectra and from the conditional variance argument 
which ensures the non-degeneracy of ${\rm Hess} (\rho_\beta )(\hat{t})$.

\end{document}